\documentclass[%
 aip,
cp,  
 amsmath,amssymb,
 preprint,%
]{revtex4-2}
\usepackage{graphicx}
\usepackage{dcolumn}
\usepackage{bm}

\usepackage[utf8]{inputenc}
\usepackage[T1]{fontenc}
\usepackage{mathptmx}

\usepackage{subfig}
\usepackage{mwe}
\usepackage[percent]{overpic}
\usepackage{amssymb,amsmath,epsfig,mathrsfs,algorithm,algorithmic,rotating}
\usepackage{mathtools,graphics,psfrag,graphicx,color,stmaryrd}
\usepackage[table]{xcolor}
\usepackage{tikz}
\usetikzlibrary{arrows,shapes,backgrounds}
\usepackage{booktabs}
\usepackage{pgfplots}
\pgfplotsset{compat=newest}
\usepackage{slashed}
\usepackage{pgfplotstable}
\usepackage{array}
\usepackage{colortbl}
\usepackage{booktabs}
\usepackage{eurosym}
\usepackage{amsmath}
\usepackage{pgfplotstable}
\usepackage{pgfplots}
\pgfplotsset{compat=newest}
\usepackage{verbatim}
\usepackage[normalem]{ulem} 
\usepackage{fancyhdr}
\usepackage{hyperref}
\hypersetup{
pagebackref=true, 
hyperindex=true, 
colorlinks=true, 
breaklinks=true, 
urlcolor=blue, 
linkcolor=blue, 
citecolor=blue, 
filecolor=blue,
pdftitle={FEM for biomembranes}, 
pdfauthor={Aymen Laadhari}, 
pdfsubject={FEM for biomembranes} 
}

\newcommand{\tr}{\operatorname*{tr}}
\newcommand\Grad{\nabla}

\newcommand\bn{\boldsymbol{n}}

\newcommand\bD{\text{\bf D}}

\newcommand\bsigma{\boldsymbol{\sigma}}

\newcommand\bx{\boldsymbol{x}}
\newcommand\bu{\boldsymbol{u}}
\newcommand\bT{\boldsymbol{T}}
\newcommand\bv{\boldsymbol{v}}

\newcommand{\Frac}[2] {\displaystyle \frac{#1}{#2}}

\newcommand{\Th}{\mathscr{T}_{h}}
\newcommand{\T}{\mathscr{T}_{h}}

\def\bId{\text{\bf Id}}
\newcommand\ddiv{\mathop{{\rm div}}\nolimits}
\def\bnu{\boldsymbol{\nu}}

\def\Ca{\text{Ca}}
\def\Re{\text{Re}}

\def\Om{\Omega}
\def\G{\Gamma}

\def\RR{\mathbb{R}}
\def\bT{\mathbf{T}}
\def\bI{\boldsymbol{I}}

\def\Om{\Omega}
\def\G{\Gamma}

\def\RR{\mathbb{R}}
\def\bT{\mathbf{T}}

\def\bI{\boldsymbol{I}}

\usepackage{marginnote}

\newcommand{\cblue}[1]{\textcolor{black}{#1}}
\newcommand{\cred}[1]{\textcolor{black}{#1}}

\colorlet{cgray}{gray!20!white}
\colorlet{cblue}{blue!20!white}
\colorlet{RoyalBlue}{blue!20!white}
\colorlet{RoyalGreen}{green!20!white}

\newcommand{\bnabla}    { \boldsymbol{\nabla} }

\renewcommand{\leq}    {\leqslant}

\everymath{\displaystyle}

\allowdisplaybreaks

\begin{document}

\preprint{PREPRINT}

\title{\cblue{A Finite Element Approach For Modeling Biomembranes In Incompressible Power-Law Flow}}

\author{Aymen Laadhari}
\email[Corresponding author: ]{Aymen.Laadhari@ku.ac.ae}
\author{Ahmad Deeb} 
\affiliation{Department of Mathematics, College of Arts and Sciences, Khalifa University of Science and Technology, Abu Dhabi, United Arab Emirates.
}

\date{December 12, 2022} 

\begin{abstract}
We present a numerical method to model the dynamics of inextensible biomembranes in a quasi-Newtonian  incompressible 
flow, which better describes hemorheology in the small vasculature.
We consider a level set model for the fluid-membrane coupling,
while the local inextensibility condition is relaxed by introducing a penalty term.
The penalty method is straightforward to implement from
any Navier-Stokes/level set solver and allows substantial computational savings over a mixed formulation.
A standard Galerkin finite element framework is used with an arbitrarily high order polynomial approximation for better accuracy in computing the bending force.
The PDE system is solved using a partitioned strongly coupled scheme based on Crank-Nicolson time integration.
Numerical experiments are provided to validate and assess the main features of the method.
\end{abstract}

\maketitle

\section{Introduction}
This paper is concerned with the numerical study of the time-dependent dynamics of biomembranes in a surrounding Newtonian and non-Newtonian flow.
The coupled fluid-membrane problem is highly nonlinear and time consuming.

Blood is a very complex fluid. Its rheology at the macroscopic scale depends both on the individual dynamics of its embedded entities and their fluid-structure interactions at the microscopic level.
\cblue{Red blood cells, referred to as RBCs,} represent its main cellular component;
They are responsible for the supply of oxygen and the capture of carbon dioxide.
In the laboratory, giant unilamellar vesicles (diameter $\approx10{\mu}$m) are biomimetic artificial liquid drops, used in vitro and in silico to study the RBCs.
Understanding the dynamics of RBCs in flow remain a difficult problem in the field of computational physics and at the theoretical level as well, consequently leading to a growing
interest in the past two decades. 
In the published literature, several works have covered the areas
of experimental biology~\cite{Song2015},
theoretical biology~\cite{Safran}, physics~\cite{badr,FLM:389603,PhysRevE.92.012717}
and applied mathematics~\cite{dziuk,Nurnberg2015}.

%
%
From a mechanical continuum perspective, Canham~\cite{Can-1970}, Helfrich~\cite{Hel-1973} and Evans~\cite{Evans1974923} independently introduced in the early 1970s a model to describe the mechanics of lipid bilayer membranes, where cellular deformations are driven by the principal curvatures.
This results in a highly nonlinear membrane force with respect to shape, see a mathematical derivation for a generalized energy functional based on shape optimization in~\cite{laadhari10}.

Different methods have been developed to study the dynamics of biomembranes in a Newtonian flow.
We can distinguish the level set method \cite{CotMaiMil2008,LSM2014,Doyeux2013251},
the phase field method~\cite{DuLiuWan-2004},
the immersed boundary method~\cite{Hu2014670}, 
the boundary integral method~\cite{RahVeeBir-2010},
the parametric finite elements~\cite{Nurnberg2015},
and the lattice Boltzmann method~\cite{Kaoui2011}.
%
From a numerical point of view, iterative and fully explicit decoupling strategies for the membrane-fluid problem are the most used techniques~\cite{Salac2012,Doyeux2013251}.
An explicit treatment of the bending force usually leads to numerical instability problems and severe time step limitations, depending on the local mesh size and bending stiffness.
However, only few works devised semi-implicit~\cite{Nurnberg2015}
or fully implicit time integration schemes~\cite{LSM2014,LSMS2018}.
Although stability is improved, a high computational burden is generally obtained with implicit strategies.
Other interesting decoupling strategies can be found in~\cite{VALIZADEH2022114191,LAADHARI2018376,LAADHARI201835,torres2019}.

While blood flow behaves like Newtonian fluid in larger diameter arteries at high shear rates, it exhibits non-Newtonian behavior in small diameter arteries with low shear rates at the microscopic scale~\cite{cokelet1963rheology}.
Non-Newtonian rheology is mainly due to polymerization and the underlying mechanisms leading to the activation and deactivation of platelets and the interactions between different microscopic entities.
Blood viscosity tends to increase at low shear rates as RBCs aggregate into a roller shape.
The Casson, Power-Law, and Quemada models are the most widely used generalised Newtonian rheologies for blood~\cite{NEOFYTOU2003127,copley1960apparent}. To our knowledge, such models have not yet been studied for the current problem.
In this work, we consider a quasi-Newtonian power law model to describe the hemorheology.


The aim of this paper is to study the dynamics of biomembranes in a complex non-Newtonian incompressible viscous flow.
In order to keep a reasonable computational cost compared to a fully mixed formulation, we design a penalty method to account for the local inextensibility of the membrane.
Various higher-order finite element approximations are used to better approximate the bending force.
We present a set of numerical examples to validate and show the main features of the method.

%
\begin{figure}[!ht]
  \centering
   \includegraphics[width=5.5cm,height=5.5cm]{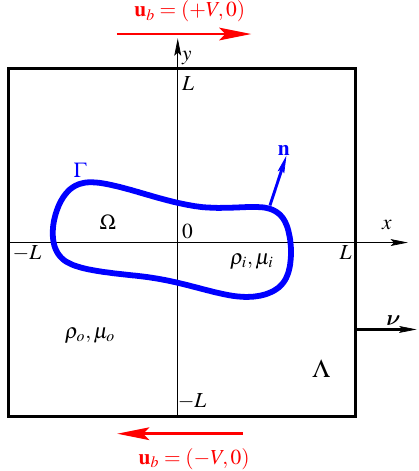}
  \caption{Sketch of the membrane $\Gamma$ embedded into a computational domain $\Lambda$, while $\Omega$ is the inner region.}
  \label{fig:dom}
\end{figure}
\section{Mathematical setting}\label{sec:model}
\subsection{Membrane model}
The deformations of the membrane allow minimizing the Canham-Helfrich-Evans \cite{Can-1970,Hel-1973} bending energy while preserving the local inextensibility of the membrane. 
Let $H$ be the mean curvature, corresponding to the sum of the principal curvatures on the membrane.
In the two-dimensional case, the membrane minimizes the bending energy given by:\begin{equation}
\mathsf{J}(\Omega)=\frac{k_b}{2}\int_{\partial\Omega}\left(H(\Omega)\right)^2 {\rm d}s,
\label{energy}\end{equation} 
where
$k_b \approx10^{-20}/10^{-19}kg\,m^2\,s^{-2}$
is the bending rigidity modulus.
The energy is a variant of the Willmore energy~\cite{Wil-1993}.
Let $T$ be the final time of the experiment.
For any time $t\in\, [0,T]$, $\Omega(t)\subset\RR^d$, $d=2,3$, is the interior domain of the membrane $\Gamma(t)=\partial\Omega(t)$, assumed Lipschitz continuous.
The membrane is embedded in the domain $\Lambda$ which is large enough so that $\G(t)\cap\partial\Lambda=\emptyset$, see Fig.~\ref{fig:dom}.
Hereafter, the dependence of $\Omega$ and $\Gamma$ upon t is dropped to alleviate notations.

For a membrane with fixed topology, 
the Gauss-Bonnet theorem~\cite{feng-klug-06} states that the energy term weighted by $k_g$ is constant and can be ignored. 
The spontaneous curvature helps describe the asymmetry of phospholipid bilayers at rest, e.g. when different chemical environments exist on either side of the membrane. We assume $H_0=0$.
Let $\bn$ and $\boldsymbol\nu$ be the  outward unit normal vector  $\Gamma(t)$ and on $\partial\Lambda$, respectively.
We introduce the surface gradient $\bnabla_s\cdot= (\bId-\bn\otimes\bn)\mbox{ }\bnabla\cdot$,
 surface divergence $\ddiv_s\cdot=\tr(\bnabla_s\cdot)$
and surface Laplacian $\Delta_s\cdot= \ddiv_s\left(\bnabla_s \cdot\right)$,
where $\bId$ is the identity tensor.
%
The expression and derivation of the bending force using shape optimization tools can be found in~\cite{laadhari10}. 

Membrane deformations are subject to specific constraints.
\cblue{Fluid incompressibility is assumed,}
this is ${\rm div}\,\bu=0\mbox{ in }\Lambda$.
In addition, RBCs are phospholipid bilayers  with local membrane inextensibility.
This corresponds to a zero surface divergence, i.e.
${\rm div}_s\bu=0$ over $\G$, that helps preserve the local perimeter.
Global perimeter conservation follows from Reynolds' lemma~\cite{LSM2014}.
As a consequence, a saddle point formulation results in a membrane surface force that balances the jump in hydrodynamic stress tensor  and appears in the right side of~\eqref{eq-pb-7}.
\subsection{Level set description}
The motion of the membrane is followed implicitly in a level set framework as the zero level set of a function~$\varphi$.
For $t\in\,]0,T[$, $\varphi$ is initialized by a signed distance $\varphi_0$ to $\G(0)$ and satisfies the transport equation~\eqref{eq-pb-1},
with $\bu$ the advection vector and $\varphi=\varphi_b$ on the upstream boundary
$\displaystyle\Sigma_{-}  = \left\lbrace\bx \in \partial\Lambda : \bu\cdot\bnu(\bx) < 0\right\rbrace$.
Geometric quantities such as $\bn=\Grad\varphi / |\Grad\varphi|$, $H=\ddiv_s\bn$ and bending force are
coded in terms of $\varphi$ and are then extended to the entire computational domain $\Lambda$.
Over time, a redistancing problem is resolved to maintain the signed distance property lost by advection~\cite{Laadhari20171047}.
Indeed, a too large or too small gradient of $\varphi$ close to $\G$ deteriorates the precise compting of the surface terms.
Let $\varepsilon$ be  a regularization parameter.
We introduce the regularized Heaviside $\mathscr{H}_\varepsilon$ and Dirac $\delta_\varepsilon$ functions:
\\
\begin{eqnarray*}
\mathscr{H}_\varepsilon(\varphi)
 = \left\{
\begin{array}{ll}
 0, &
 \mbox{ when } \varphi<-\varepsilon\\
 \Frac{1}{2}\left( 1 + \Frac{\varphi} {\varepsilon}
 +\Frac{1}{\pi} \sin\left(\Frac{\pi\varphi}{\varepsilon}\right)\right), & \mbox{ when } \left|\varphi\right| \leq \varepsilon,\\
 1, & \mbox{ otherwise}
 \end{array} \right.
 &\quad
\text{ and }
\quad \delta_{\varepsilon}(\varphi) 
 = \Frac{{\rm d}\mathscr{H}_\varepsilon}{{\rm d}\varphi} (\varphi).
\end{eqnarray*}
Given a function $\zeta$ defined on $\G$ and its extension $\tilde{\zeta}$ to $\Lambda$, surface integrals are approximated as follows:
\begin{equation*}\label{eq:GL}
\int_{\G} \zeta(\bx)\,{\rm d}s 
\approx
\int_{\Lambda}|\nabla \varphi|\,\delta_\varepsilon\left(\varphi\right)\,\tilde{\zeta}(\bx)\,{\rm d}x.
\end{equation*}
%
%
\subsection{Governing equations}
We assume constant densities $\rho_{i}$ and $\rho_{o}$ 
 inside and outside of the membrane, respectively.
%
Let us introduce the fluid velocity $\bu$ and the pressure $p$ which represent a Lagrange multiplier corresponding to the incompressibility constraint on $\Lambda$.
Analogously, a position-dependent surface tension $\lambda$ helps imposing the local inextensibility constraint on $\G$.
Let $\bD(\bu)=(\bnabla\bu+\bnabla\bu^T)/2$ be the shear strain rate tensor, 
so the fluid Cauchy stress tensor is 
$\bsigma = \bT-p\bI$ where $\bT$ is the stress deviator.
%
The normal stress jump $[\bsigma\bn]_-^+=\bsigma_+\bn-\bsigma_-\bn$ on $\G$ describes the interactions of the membrane with the surrounding fluid~\cite{LSMS2018}, while the stress discontinuity is calibrated by~\eqref{eq-pb-5}.
For a simple shear flow, $\bu_b$ is the shear rate on $\Sigma_D\subset\partial\Lambda$, while natural boundary conditions are prescribed on $\Sigma_N\subset\partial\Lambda$.

We assume a quasi-Newtonian power-law model~\cite{NEOFYTOU2003127}
where the nonlinear constitutive equation expresses the stress deviator with a power-law viscosity function as
\begin{equation}\label{eq:powerlaw}
\bT = 2  \eta\left( \left| \bD({\bu}) \right|^2 \right)  \bD({\bu}), \mbox{ with  } \eta \left( \gamma \right) = 
K \gamma^{\displaystyle (\upsilon-1)/2},\,
 \mbox{ for all  }
 \gamma\in\mathbb{R},
 \end{equation}
%
where $\upsilon>0$ and $K$ are the power index and consistency index, respectively.
According to~\cite{Walburn1976}, $\upsilon=0.7755<1$ (i.e. a shear thinning fluid) and $K=14.67\times 10^{-3}$ Pa s  for normal blood samples obtained using a multiple regression technique.
The Newtonian case $\upsilon=1$ corresponds to a linear stress-strain relationship that 
reduces the viscosity function  $\eta( \gamma)=K$ to a constant. 
By analogy with the Newtonian case, $K=\mu_i$ and $K=\mu_o$ stand for the values of the consistency index in the intra- and extra-membrane domains, respectively.
%
%

We perform a dimensionless analysis.
Let U be the maximum velocity on $\Sigma_D$ and $D$ the diameter of a circle having the same membrane perimeter.
We consider the dimensionless Reynolds number $\Re=\rho_o U D\mu_o^{-1}$ which expresses the ratio between the inertial and viscous forces, and the capillary number $\Ca=\mu_oD^2Uk_b^{-1 } $ which compares the flow force to the bending resistance of the membrane.
Furthermore, the 
parameter $\beta=\mu_i/\mu_o$ represents the ratio of consistency indices and corresponds to the viscosity ratio with respect to extracellular viscosity in the Newtonian case.
The regularized dimensionless viscosity function is:
\[
    \mu_\varepsilon(\varphi)  | \bD({\bu}) |^{\upsilon-1} 
    = \left(
    \mathscr{H}_\varepsilon(\varphi) + \beta \left(1-\mathscr{H}_\varepsilon(\varphi)\right)
    \right)   | \bD({\bu}) |^{\upsilon-1}.\]
Following~\cite{LSMS2018}, we choose $\rho_i=\rho_o$.
\cblue{Let $\bsigma_\varepsilon$ stand for the regularized Cauchy stress tensor}.
The dimensionless reduced area
 $\Xi_{2d}=4\pi |\Omega| / |\G|^{2}\in]0,1]$
compares the area of the interior domain to that of a circle with the same perimeter.
The dimensionless coupled problem writes: find $\varphi$, ${\bf u}$, $p$ and $\lambda$ such that
\begin{subequations}
\begin{eqnarray}
   \partial_t \varphi
    +
    {\bf u}.\nabla \varphi
    &=&
    0
    \ \mbox{ in } ]0,T[ \times \Lambda
    \label{eq-pb-1}
    \\
    \Re \,
    \left(
     \partial_t \bu
      +
      {\bf u}.\nabla {\bf u}
    \right)
    -
    {\bf div} \big(
        \cblue{\bsigma_\varepsilon}( D({\bf u}), p, \varphi)
    \big)
    &=&
    0
    \ \mbox{ in } ]0,T[ \times (\Lambda \backslash \partial\Omega)
    \label{eq-pb-2}
    \\
    {\rm div}\,{\bf u}
    &=&
    0
    \ \mbox{ in } ]0,T[ \times \Lambda
    \label{eq-pb-3}
    \\
    {\rm div}_s\,{\bf u}
    &=&
    0
    \ \mbox{ on } ]0,T[ \times \partial\Omega
    \label{eq-pb-6}
    \\
    \left[{\bf u}\right]_-^+ 
    &=&
    0
    \ \mbox{ on } ]0,T[ \times \partial\Omega
    \label{eq-pb-4}
    \\
    \left[\bsigma_\varepsilon\bn\right]_-^+ 
    &=&
   \bnabla_s \lambda - \lambda H \bn +
   (2\Ca)^{-1} \left(2\Delta_sH+H^3\right)\bn
    \ \mbox{ on } ]0,T[ \times \partial\Omega
    \label{eq-pb-5}
    \\
    \varphi &=& \varphi_b
	\ \mbox{ on }\ 
        ]0,T[\times \Sigma_{-}
    \label{eq-pb-7}
    \\
    {\bf u} &=& {\bf u}_b
	\ \mbox{ on }\ 
        ]0,T[\times \Sigma_D
    \label{eq-pb-8}
    \\
   \bsigma.\boldsymbol{\nu} &=& 0
	\ \ \mbox{ on }\ 
        ]0,T[\times \Sigma_N
    \label{eq-pb-9}
    \\
    \varphi(0) &=& \varphi_0
	\ \mbox{ in }\ 
        \Lambda
    \label{eq-pb-10}
    \\
    {\bf u}(0) &=& {\bf u}_0
	\ \mbox{ in }\ 
        \Lambda.
    \label{eq-pb-11}
\end{eqnarray}
\end{subequations}
%
Let  ${\epsilon_\lambda}=10^{-8}$ be the penaly parameter.
To make the method straightforward to implement from any Level Set / Navier-Stokes solver
and considerably reduce the size of the linear system to be solved,
the inextensibility constraint is relaxed by introducing a penalty term.
Indeed, the corresponding minimization problem should be approximated by another minimization
problem by penalizing the local inextensibility constraint for the velocity~\eqref{eq-pb-6}. 
See analogous penalty method for other applications in~\cite{janela}.

To overcome instability problems when solving the level set equation using the standard Galerkin method, there are a variety of stabilization methods such as the streamline diffusion method, the subgrid viscosity method and the Streamline Upwind Petrov-Galerkin (SUPG) method used in this work.
The latter introduces a stabilization term by adding a diffusion in the streamline direction.

We introduce the functional spaces of admissible velocity $\bu$, pressure $p$ and level set $\varphi$:
\begin{align*}
 &\mathbb{V}(\bu_b)=
\Big\{
  \bv \in \left( H^1\left(\Lambda\right)\right)^d : 
  \bv = \bu_b, \text{ on } \Sigma_D
  \Big\},
  \qquad
  \mathbb{Q}=  \left\{
  q \in  L^2\left(\Lambda\right) :  
  \int_\Om q  =0 \right\},
  \\
  &
   \mathbb{X}(\varphi_b)
 =  \left\{
 \psi \in  W^{1,\infty}\left(\Lambda\right)\cap H^1\left(\Lambda\right)\, : \ \psi = \varphi_b, \text{ on } \Sigma_-\right\}.
\end{align*}

To reduce a derivation order of $\varphi$ when evaluating the bending strength, we
use the Green formula on a closed surface. See e.g.~\cite{LSM2014}.
%
%
Testing with appropriate test functions and integrating~\eqref{eq-pb-2} over $\Omega$ and $\Lambda\backslash\overline{\Omega}$ separately, the variational problem writes:\\
    Find
    $\bu \in \mathcal{C}^0\Big( ]0,T[, L^2(\Lambda)^d\Big)\cap L^2\Big( ]0,T[,\mathbb{V}(\bu_b) \Big)$, 
    $p \in L^2\Big( ]0,T[,\mathbb{Q} \Big)$,
    and
    $\varphi \in \mathcal{C}^0\Big( ]0,T[, L^2(\Lambda)^d\Big) \cap L^2\Big( ]0,T[,\mathbb{X}\left(\varphi_b\right) \Big)$ 
    such that
\begin{subequations}
\begin{eqnarray}
      &&\Re 
       \int_{\Lambda} \left( \frac{\partial\bu}{dt} + \bu\cdot\Grad\bu \right)\cdot   \bv   +
       \int_{\Lambda} 2\mu_\varepsilon(\varphi)  | \bD({\bu}) |^{\upsilon-1}  \bD({\bu}) : \bD(\bv)
       +\frac{1}{\epsilon_\lambda} \int_{\Lambda}  \ddiv_s(\bu) \, \ddiv_s(\bv)
                                        |\Grad\varphi|  \delta_{\varepsilon}(\varphi) 
       \nonumber\\
     && \quad
        -\int_{\Lambda} p \, \ddiv \, \bv
        +\frac{1}{2\Ca}\int_\Lambda 
         \delta_{\varepsilon}(\varphi) |\Grad\varphi|
        \Big( 
                                    2\bnabla_{s} H \cdot\bnabla_{s} (\bn\cdot\bv) 
                                    - H^3 \bn\cdot\bv
                                         \Big)
       =
       \int_{\Sigma_N}\bsigma\bnu\cdot\bv, 
       \hspace{1.5cm} \forall \bv \in \mathbb{V}(0),\label{eq:vf1}\\
     && \int_{\Lambda} q\ \ddiv\, \bu = 0 , 
       \hspace{11.25cm}  \forall q \in \mathbb{Q},
       \label{eq:vf2}
       \\
     && \int_{\Lambda} \frac{\partial\varphi}{\partial t}\psi + \int_{\Lambda} \left(\bu\cdot\nabla\varphi\right)\psi 
         +  \int_{\Lambda} \xi \left( \tau; \varphi, \psi \right)
= 0,
     \hspace{6.45cm} \forall \psi \in \mathbb{X}\left(0\right). \label{eq:vf4}
\end{eqnarray}
\end{subequations}
Here, $ \xi \left( \tau; \varphi, \psi \right)$ stands for the SUPG stabilisation term and 
$\tau$ is a stabilization parameter defined element wise to control the amount of diffusion.
\section{Numerical approach}\label{sec:num}
The interval $[0,T]$ is divided into $N$ sub-intervals $[t^n,t^{n+1})$ with $0\leq n \leq N-1$ of constant step $\Delta t$ .
For $n>0$, $\bu^{n}$, $p^{n}$ and $\varphi^{n}$ are computed by induction to approximate
 $\bu$, $\varphi$ and $p$ at $t^n$.
We use the Crank-Nicolson scheme for the time discretization of \eqref{eq-pb-1} and \eqref{eq-pb-2} without the need to bootstrap the initial conditions. \cred{The choice of this scheme was for its simplicity to implement and being a second order one-step integrator.}
The discretized~\eqref{eq-pb-1} writes
\[ \varphi^{n+1} = \varphi^n + \frac{\Delta t}{2}\left(
             {\bf u}\cdot\nabla \varphi^{n+1} +
             {\bf u}\cdot\nabla \varphi^{n}
             \right)
    \quad  \ \mbox{ in } \Lambda.
\]
For the spatial discretization, we consider a partition $\Th$ of $\Lambda$ consisting of geometrically conformal
open simplicial elements $K$.
We define the mesh size as the diameter of the largest mesh element $h=\max h_K$ with $K\in\Th$.

We consider a Taylor-Hood finite element approximation for $\bu$ and $p$.
After using a surface Green's transformation, the evaluation of the Canham-Helfrich-Evans force requires a third-order derivative in $\varphi$ which induces numerical oscillations when using lower-order polynomial approximations.
 To avoid introducing additional mixed variables and additional equations as in \cite{LSM2014}, higher degree polynomials are considered for the discretization of $\varphi$ because the bending force requires its fourth order derivatives.
For the SUPG method, the streamline diffusion parameter is chosen numerically proportional to the local mesh size, this is
$\tau_K = C h_K / \max\left\{|\bu|_{0,\infty,K}, \text{tol}/h_K\right\}$,
where C is a scaling constant and $\text{tol}/h_K$ helps to avoid division by zero.
To overcome the instability problems induced by an explicit decoupling, we consider a partitioned implicit strategy based on a fixed point algorithm, as detailed in Alg.~\ref{alg:1}.
\begin{algorithm}[H]
{\small  \caption{Fluid-membrane coupling}
\label{alg:1}
 \begin{algorithmic}[1]
 \STATE $n=0$:
    let $\varphi^{0}$ and ${\bf u}^{0}$ being given
     \FOR{$n=0,\dots,N-1$}
     \STATE  Initialize $\bu^{n+1,0}=\bu^{n}$, $\varphi^{n+1,0}=\varphi^{n}$
 \WHILE{$e^k< 10^{-6}$ }
     \STATE Compute $\varphi^{n+1,k+1}$ using  $\bu^{n+1,k}$
     \STATE Compute $\bu^{n+1,k+1}$, $p^{n+1,k+1}$ using  $\varphi^{n+1,k+1}$
     \STATE Compute the error $e^k= |\bu^{n+1,k+1}-\bu^{n+1,k}|_{1,2,\Lambda}/|\bu^{n+1,k}|_{0,2,\Lambda}
      + 
       |\varphi^{n+1,k+1}-\varphi^{n+1,k}|_{0,2,\Lambda}/|\varphi^{n+1,k}|_{0,2,\Lambda}$
 \ENDWHILE
     \STATE  Update $\bu^{n+1}=\bu^{n+1,k+1}$, $\varphi^{n+1}=\varphi^{n+1,k+1}$
 \ENDFOR
 \end{algorithmic}}
 \end{algorithm}
%
%
%
\section{Numerical examples}\label{sec:sim}
%
\begin{figure}[ht]
    \centering
    \begin{minipage}{0.65\textwidth}
     \includegraphics[width=0.32\textwidth]{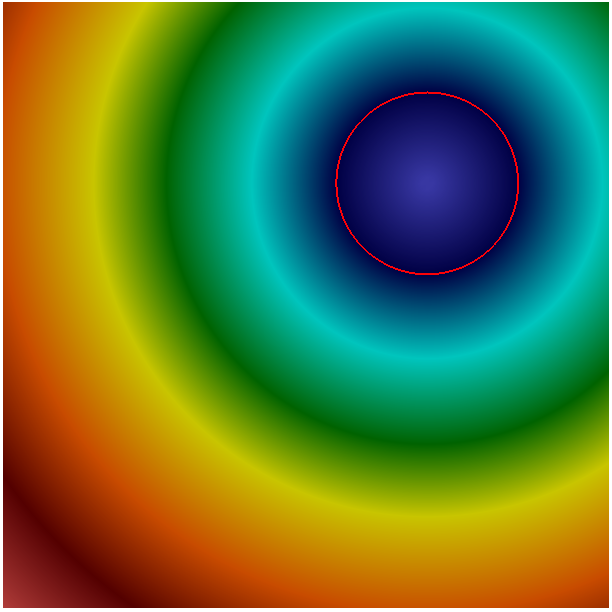}
     \hspace{-0.cm}
     \includegraphics[width=0.32\textwidth]{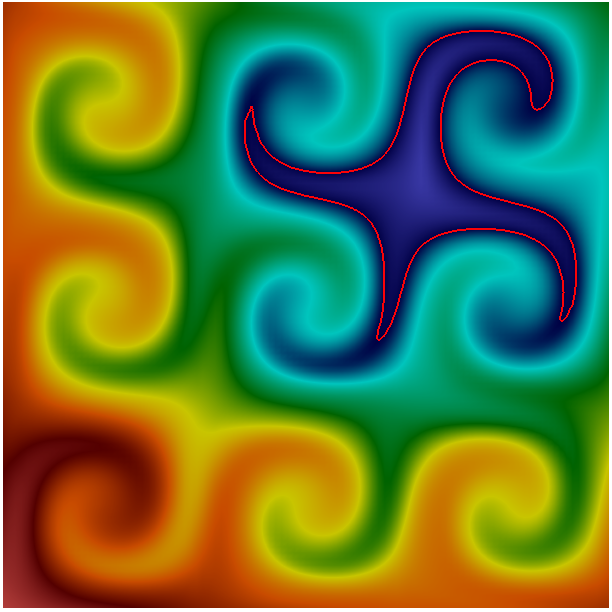}
     \hspace{-0.cm}
     \includegraphics[width=0.32\textwidth]{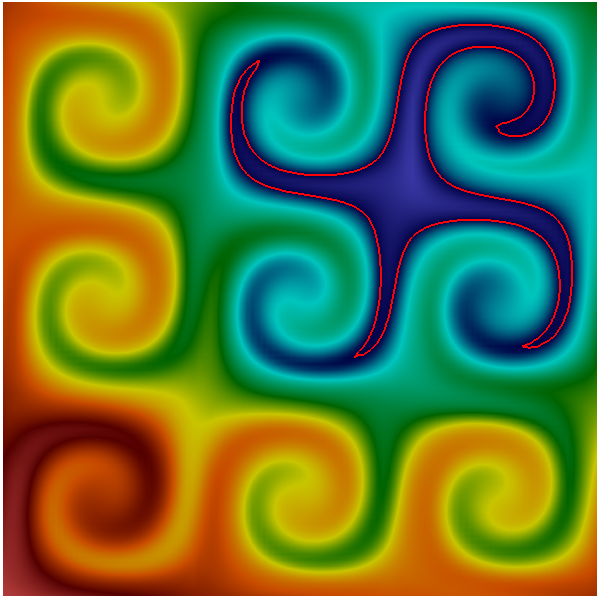}
     \newline
     \includegraphics[width=0.32\textwidth]{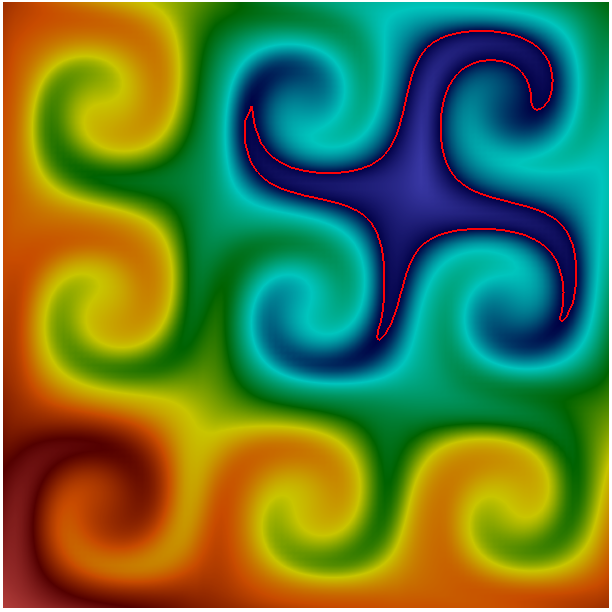}
     \hspace{-0.cm}
     \includegraphics[width=0.32\textwidth]{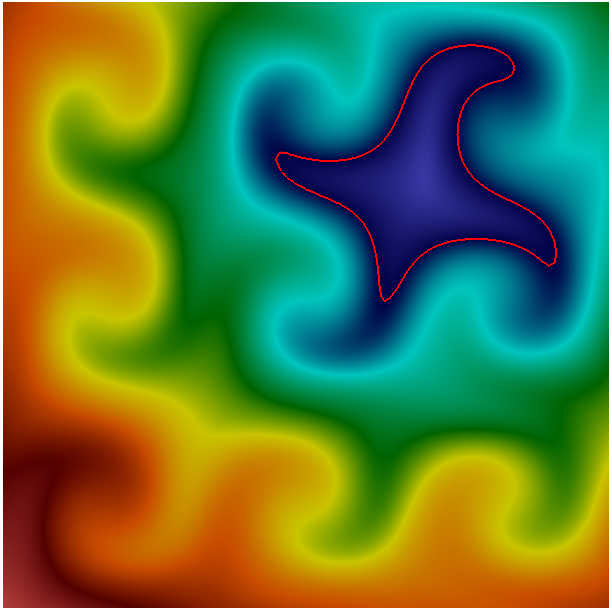}
     \hspace{-0.cm}
     \includegraphics[width=0.32\textwidth]{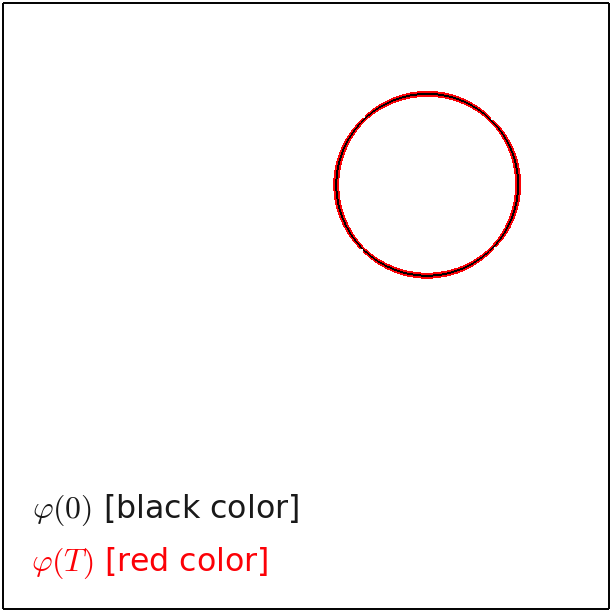}
    \end{minipage}
    \hfill
    \begin{minipage}{0.34\textwidth}
    \hspace{-0.1cm}
     \includegraphics[width=0.96\textwidth]{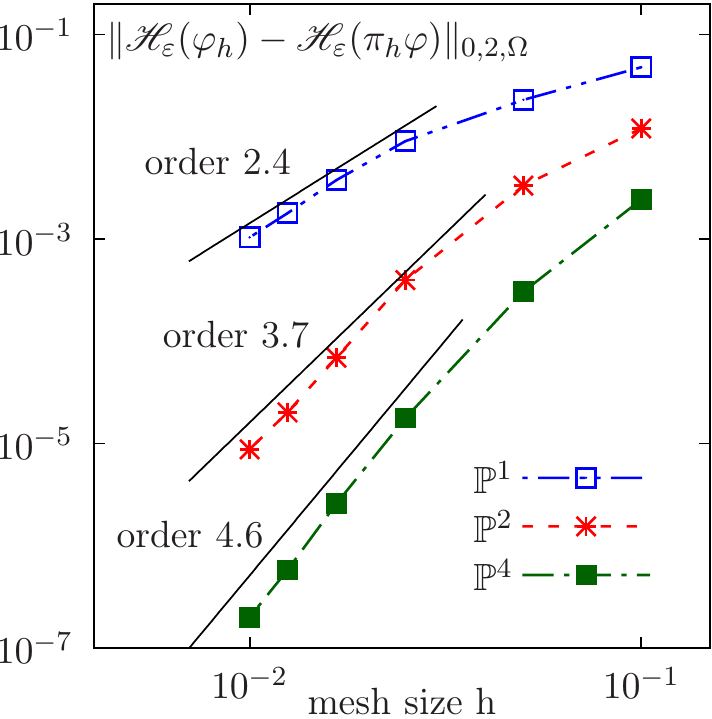}
     \newline
     \centering
      $ $\newline
     \includegraphics[height=1.0cm]{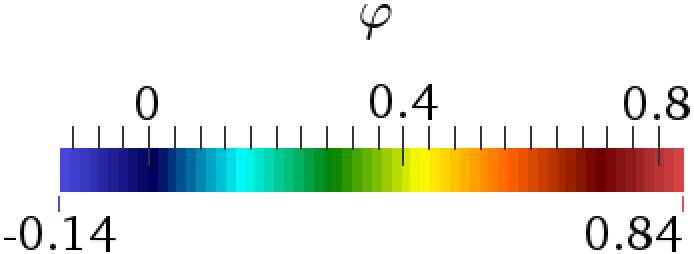}
    \end{minipage}
\caption{Reversible vortex. (Left) Snapshots showing the interface deformations at 
$t\in\{0, 0.25, 0.57,$ $0.75, 0.875, 1\}$ with $h=0.01$. 
(Right) Spatial convergence in $L^2$ norm for high-order 
finite element approximations.}
\label{fig:test}
\end{figure}
\subsection{Example 1: Reversible Vortex - Grid convergence.}
Simulations were performed using FEniCSx\cite{AlnaesEtal2015}. 
%
To evaluate the capability of the level set solver for high-order finite elements, 
necessary afterwards for an accurate assessment of highly nonlinear bending force, 
we consider a reversible vortex test case featuring large deformations of the interface. 
The computational domain is $\Lambda=[0,1]^2$.
A circular interface of radius $R=0.15$ initially centered at $(0.7,0.7)$ is stretched into thin
filaments which are coiled like a starfish by a vortex flow field.
The deformations are periodic and the stretching of the membrane 
unravels before the interface regains its circular shape after a period at $t=T$. 
The maximal deformation $\psi$ occur at $t=T/2$, with $\psi=3$ and $T=1$ in numerical computations.
Similar 2D and 3D test cases are widely used to test interface tracking methods.
We follow LeVeque's test \cite[Example 9.5]{Leveque-1996} and consider a velocity field at $\bx=(x,y)^T\in\Lambda$ given by
$$\bu(t,\bx) =\left( 
-2 \sin(\psi \pi x)^2  \ sin(\psi \pi y) \cos(\psi \pi y) \cos(\pi t/T)
,
 2 \sin(\psi \pi y)^2   \sin(\psi \pi x) \cos(\psi \pi x) \cos(\pi t/T)
 \right)^T.$$


The spatial accuracy of the finite element numerical approximations is studied by computing the errors in $L^2(\Lambda)$
norm on successively refined meshes with respect to an exact reference solution $\pi_h\varphi$ at $t=T$,
where $\pi_h$ represents the Lagrange interpolation operator.
Errors are calculated after one stretching period.
For $k$ the degree of the polynomial approximation, the time step $\Delta t=h^k$ is chosen small enough not to significantly influence the overall accuracy.
Fig. \ref{fig:test} reports the convergence of calculated errors with respect to the mesh size for several polynomial
finite element approximations.
Convergence rates are also displayed, showing for instance an almost second-order accuracy for $k=1$
and fifth-order accuracy for $k=4$.
%
\begin{figure}[!h]
 \centerline{
     \includegraphics[height=4.9cm]{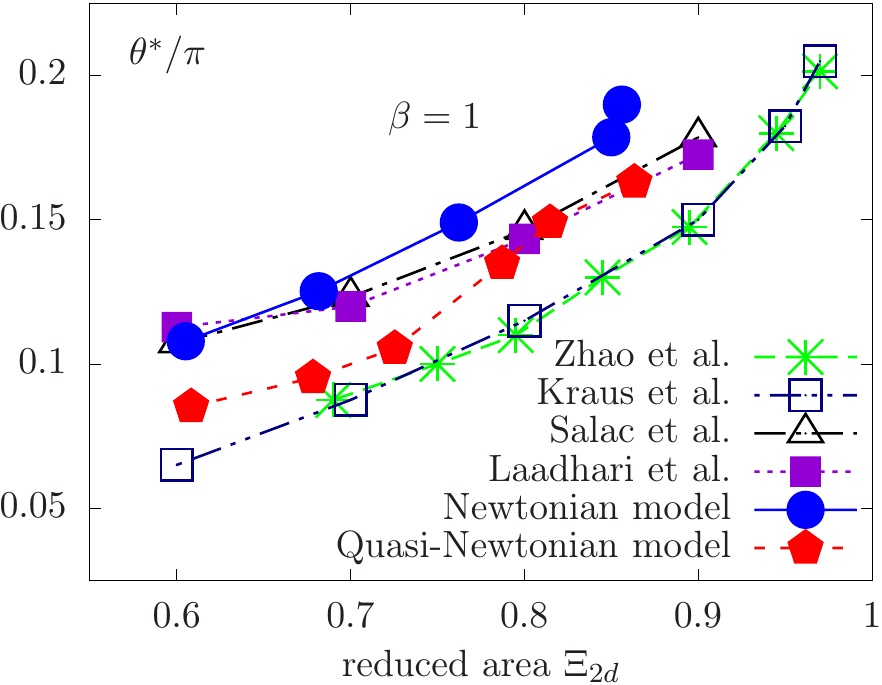}
     }
 \caption{TT regime: Change in $\theta^*/\pi$ with respect to $\Xi_{2d}$ for
    a viscosity ratio $\beta=1$.
    Comparisons  with results from~\cite{Salac2012} and~\cite{LSMS2018} ($\Re=10^{-3}$, $\Ca=100$),
     \cite{Zhao2011}
    ($\Ca=9$) and \cite{Kraus1996} ($\Ca=10$).}
 \label{fig:equiAngle}
 \centering
     \includegraphics[height=4.9cm]{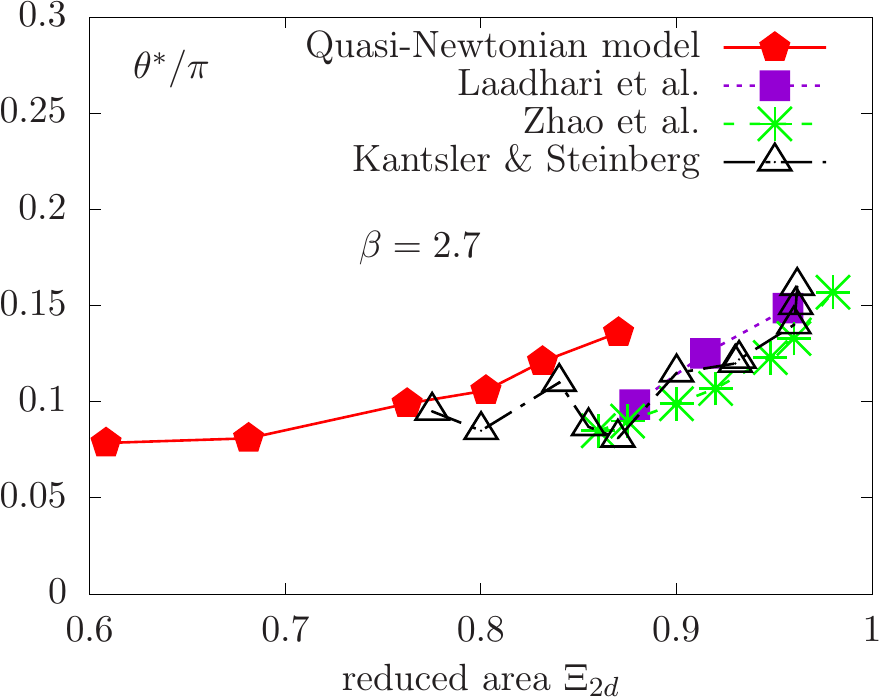}
 \caption{TT regime: Change in $\theta^*/\pi$ with respect to $\Xi_{2d}$ for 
 $\beta=2.7$. 
    Comparison of non-Newtonian model with results in~\cite{Zhao2011,LSMS2018}
    and measurements in~\cite{Kantsler2005}.}
 \label{fig:equiAngle2_7}
\end{figure}

\subsection{Example 2: Dynamics of the biomembrane in Newtonian and quasi-Newtonian flows.}

We first proceed to a quantitative validation with some experimental and numerical results available in the literature in the case of a purely Newtonian flow.
We set $\upsilon=1$, a viscosity contrast $\beta=1$, $\Ca=10^2$ and $\Re=9\times10^{-3}$.
More details on the physiological values of the Reynolds number at the level of RBCs are available in~\cite{Salac2012}.
The membrane follows a tank-treading type movement, called TT, where it reaches a steady state characterized by a fixed angle of inclination; The surrounding fluid continues its rotation tangentially to the membrane.
We consider different values of the reduced areas $\Xi_{2d}\in [0.6,1]$, and calculate the angle of inclination at equilibrium $\theta^\star$.
Fig.~\ref{fig:equiAngle} and Fig.~\ref{fig:equiAngle2_7} plot the change in $\theta^\star/\pi$ against $\Xi_{2d}$ in both Newtonian and quasi-Newtonian cases for different values of the viscosity ratio $\beta$.
The results are compared with those of Kraus et al.~\cite{Kraus1996}, Zhao et al.~\cite{Zhao2011}, Salac et al.~\cite{Salac2012} and Laadhari et al.~\cite{LSMS2018}, showing good overall consistency.
However, note that the values obtained with $\upsilon=0.7755$ fit slightly better compared to those of the Newtonian model which are a little higher than the other curves.
This cannot be confirmed at all, given the setting of different dimensionless values such as $\Re$ and $\Ca$ in the different experiments.
An in-depth study is in progress and will be the subject of a forthcoming work.
%
\begin{figure}[!h]
   \centering \includegraphics[width=15cm]{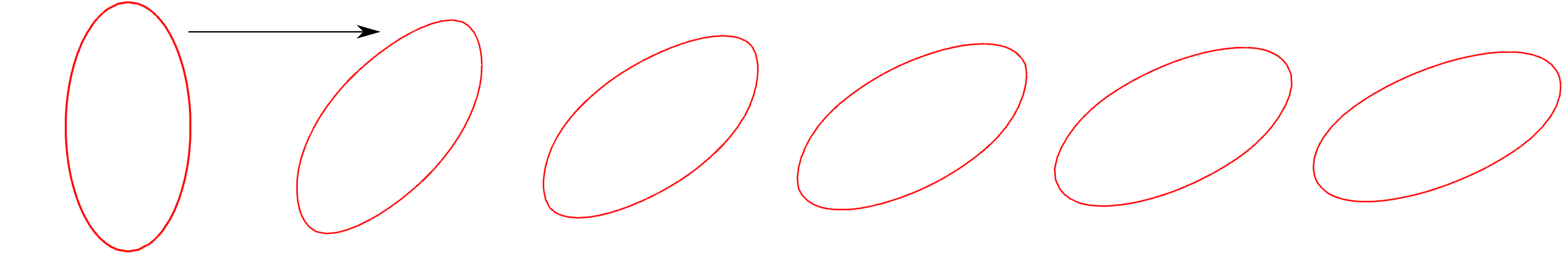}
   \caption{TT regime for $\Xi_{2d}=0.68$, $\beta=1$, $\Ca=4\times 10^4$ and $\Re=9\times 10^{ -3}$ at $t\in\big\lbrace 0, 0.125, 0.25, 0.5, 1, 2\big\rbrace$.}
   \label{TT}
\end{figure}

Simulations are now performed using a different ratio $\beta^\star=2.7$ in the non-Newtonian case.
We calculate the angle of inclination $\theta^\star/\pi$ and compare with some numerical~\cite{Zhao2011}
and experimental~\cite{Kantsler2005}
 results available only for larger reduced areas.
Fig.~\ref{fig:equiAngle}(right) shows close but slightly higher equilibrium angles when the shape of the membrane becomes close to a circle.
The deviations can be mainly due to the non-Newtonian model, but also to the different values of the confinement levels and the boundary conditions used in the different works.

\begin{figure}[!h]
   \centering
   \hspace{-0.65cm} \includegraphics[width=15cm]{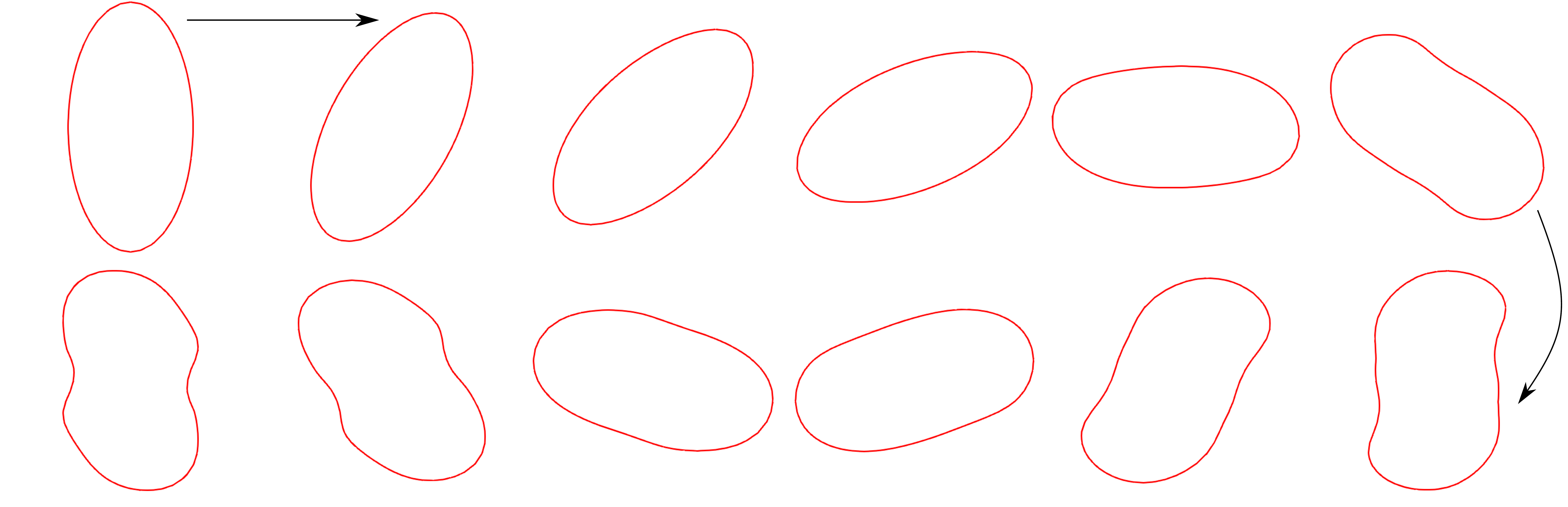}     
   \caption{Snapshots showing a membrane in TB regime for $\Xi_{2d}=0.68$, $\beta=10$, $\Ca=4\times 10^4$ and $\Re=9\times      
    10^{ -3}$, at times $t\in\big\lbrace 0, 0.13, 0.25, 0.5, 1.25, 2, 2.25, 2.38, 3,  4, 4.5, 4.7
    \big\rbrace$, respectively.
   }
   \label{TB}
\end{figure}

\cblue{According to an experimental systematic study on individual  individual red cells in a simple shear flow,
a change in dynamics occurs when the viscosity ratio exceeds a critical value depending on the reduced area~\cite{Fischer-1978}. This is the tumbling regime, noted by TB, which is characterized by the periodic rotation of the membrane around its axis.}
A well-known empirical model was developed by Keller and Skalak~\cite{FLM:389603}.
This dynamics was obtained in the simulations with the non-Newtonian model, see Fig. \ref{TT} and Fig. \ref{TB} for the snapshots of the TT and TB dynamics obtained with the same set of parameters but with $\beta=1$ and $\beta=10$, respectively.

\section{Conclusion}

We have presented in this paper a relatively simple method for simulating the dynamics of an individual red blood cell, or inextensible biological membrane in general, in a surrounding incompressible non-Newtonian flow that better describes the hemorheology in small capillaries.

We validated our framework using high-order finite element approximations in the case of a membrane in a simple shear flow.
Simulations have shown that the method is capable of capturing the basic cellular dynamics, namely the well-known tank treading and tumbling motions. This is part of a larger ongoing work to explore the dynamics of red blood cells in small capillaries, while accounting for cell elasticity~\cite{Gizzi2015157} in non-Newtonian surrounding flow.

\begin{acknowledgments}
The authors acknowledge financial support from KUST through the grant FSU-2021-027.
\end{acknowledgments}

\bibliography{mybibfile.bib}

\begin{thebibliography}{39}%
\makeatletter
\providecommand \@ifxundefined [1]{%
 \@ifx{#1\undefined}
}%
\providecommand \@ifnum [1]{%
 \ifnum #1\expandafter \@firstoftwo
 \else \expandafter \@secondoftwo
 \fi
}%
\providecommand \@ifx [1]{%
 \ifx #1\expandafter \@firstoftwo
 \else \expandafter \@secondoftwo
 \fi
}%
\providecommand \natexlab [1]{#1}%
\providecommand \enquote  [1]{``#1''}%
\providecommand \bibnamefont  [1]{#1}%
\providecommand \bibfnamefont [1]{#1}%
\providecommand \citenamefont [1]{#1}%
\providecommand \href@noop [0]{\@secondoftwo}%
\providecommand \href [0]{\begingroup \@sanitize@url \@href}%
\providecommand \@href[1]{\@@startlink{#1}\@@href}%
\providecommand \@@href[1]{\endgroup#1\@@endlink}%
\providecommand \@sanitize@url [0]{\catcode `\\12\catcode `\$12\catcode
  `\&12\catcode `\#12\catcode `\^12\catcode `\_12\catcode `\%12\relax}%
\providecommand \@@startlink[1]{}%
\providecommand \@@endlink[0]{}%
\providecommand \url  [0]{\begingroup\@sanitize@url \@url }%
\providecommand \@url [1]{\endgroup\@href {#1}{\urlprefix }}%
\providecommand \urlprefix  [0]{URL }%
\providecommand \Eprint [0]{\href }%
\providecommand \doibase [0]{http://dx.doi.org/}%
\providecommand \selectlanguage [0]{\@gobble}%
\providecommand \bibinfo  [0]{\@secondoftwo}%
\providecommand \bibfield  [0]{\@secondoftwo}%
\providecommand \translation [1]{[#1]}%
\providecommand \BibitemOpen [0]{}%
\providecommand \bibitemStop [0]{}%
\providecommand \bibitemNoStop [0]{.\EOS\space}%
\providecommand \EOS [0]{\spacefactor3000\relax}%
\providecommand \BibitemShut  [1]{\csname bibitem#1\endcsname}%
\let\auto@bib@innerbib\@empty
\bibitem [{\citenamefont {Song}\ and\ \citenamefont
  {Hashash}(2015)}]{Song2015}%
  \BibitemOpen
  \bibfield  {author} {\bibinfo {author} {\bibfnamefont {H.}~\bibnamefont
  {Song}}\ and\ \bibinfo {author} {\bibfnamefont {Y.}~\bibnamefont {Hashash}},\
  }\bibfield  {title} {\enquote {\bibinfo {title} {Characterization of
  stress-strain behaviour of red blood cells ({RBC}s), part {II}: response of
  malaria-infected {RBC}s},}\ }\href {\doibase 10.1080/17415977.2014.922072}
  {\bibfield  {journal} {\bibinfo  {journal} {Inverse Probl. Sci. Eng.}\
  }\textbf {\bibinfo {volume} {23}},\ \bibinfo {pages} {413} (\bibinfo {year}
  {2015})}\BibitemShut {NoStop}%
\bibitem [{\citenamefont {Safran}(1994)}]{Safran}%
  \BibitemOpen
  \bibfield  {author} {\bibinfo {author} {\bibfnamefont {S.~A.}\ \bibnamefont
  {Safran}},\ }\href@noop {} {\emph {\bibinfo {title} {Statistical
  Thermodynamics of Surfaces, Interfaces and Membranes}}},\ Frontier in Physics
  Vol 90\ (\bibinfo  {publisher} {Addison-Wesley Publishing Company, Reading,
  Massachusetts},\ \bibinfo {year} {1994})\BibitemShut {NoStop}%
\bibitem [{\citenamefont {Kaoui}\ \emph {et~al.}(2008)\citenamefont {Kaoui},
  \citenamefont {Ristow}, \citenamefont {Cantat}, \citenamefont {Misbah},\ and\
  \citenamefont {Zimmermann}}]{badr}%
  \BibitemOpen
  \bibfield  {author} {\bibinfo {author} {\bibfnamefont {B.}~\bibnamefont
  {Kaoui}}, \bibinfo {author} {\bibfnamefont {G.~H.}\ \bibnamefont {Ristow}},
  \bibinfo {author} {\bibfnamefont {I.}~\bibnamefont {Cantat}}, \bibinfo
  {author} {\bibfnamefont {C.}~\bibnamefont {Misbah}}, \ and\ \bibinfo {author}
  {\bibfnamefont {W.}~\bibnamefont {Zimmermann}},\ }\bibfield  {title}
  {\enquote {\bibinfo {title} {Lateral migration of a two-dimensional vesicle
  in unbounded poiseuille flow},}\ }\href@noop {} {\bibfield  {journal}
  {\bibinfo  {journal} {Phys. Rev. E}\ }\textbf {\bibinfo {volume} {77}},\
  \bibinfo {pages} {021903} (\bibinfo {year} {2008})}\BibitemShut {NoStop}%
\bibitem [{\citenamefont {S.R.Keller}\ and\ \citenamefont
  {Skalak}(1982)}]{FLM:389603}%
  \BibitemOpen
  \bibfield  {author} {\bibinfo {author} {\bibnamefont {S.R.Keller}}\ and\
  \bibinfo {author} {\bibfnamefont {R.}~\bibnamefont {Skalak}},\ }\bibfield
  {title} {\enquote {\bibinfo {title} {Motion of a tank-treading ellipsoidal
  particle in a shear flow},}\ }\href {\doibase 10.1017/S0022112082002651}
  {\bibfield  {journal} {\bibinfo  {journal} {J. Fluid Mech.}\ }\textbf
  {\bibinfo {volume} {120}},\ \bibinfo {pages} {27} (\bibinfo {year}
  {1982})}\BibitemShut {NoStop}%
\bibitem [{\citenamefont {Choi}\ \emph {et~al.}(2015)\citenamefont {Choi},
  \citenamefont {Yi},\ and\ \citenamefont {Kim}}]{PhysRevE.92.012717}%
  \BibitemOpen
  \bibfield  {author} {\bibinfo {author} {\bibfnamefont {W.}~\bibnamefont
  {Choi}}, \bibinfo {author} {\bibfnamefont {J.}~\bibnamefont {Yi}}, \ and\
  \bibinfo {author} {\bibfnamefont {Y.}~\bibnamefont {Kim}},\ }\bibfield
  {title} {\enquote {\bibinfo {title} {Fluctuations of red blood cell
  membranes: The role of the cytoskeleton},}\ }\href {\doibase
  10.1103/PhysRevE.92.012717} {\bibfield  {journal} {\bibinfo  {journal} {Phys.
  Rev. E}\ }\textbf {\bibinfo {volume} {92}},\ \bibinfo {pages} {012717}
  (\bibinfo {year} {2015})}\BibitemShut {NoStop}%
\bibitem [{\citenamefont {Dziuk}(2008)}]{dziuk}%
  \BibitemOpen
  \bibfield  {author} {\bibinfo {author} {\bibfnamefont {G.}~\bibnamefont
  {Dziuk}},\ }\bibfield  {title} {\enquote {\bibinfo {title} {Computational
  parametric {W}illmore flow},}\ }\href {\doibase 10.1007/s00211-008-0179-1}
  {\bibfield  {journal} {\bibinfo  {journal} {Numer. Math.}\ }\textbf {\bibinfo
  {volume} {111}},\ \bibinfo {pages} {55} (\bibinfo {year} {2008})}\BibitemShut
  {NoStop}%
\bibitem [{\citenamefont {Barrett}\ \emph {et~al.}(2015)\citenamefont
  {Barrett}, \citenamefont {Garcke},\ and\ \citenamefont
  {Nurnberg}}]{Nurnberg2015}%
  \BibitemOpen
  \bibfield  {author} {\bibinfo {author} {\bibfnamefont {J.~W.}\ \bibnamefont
  {Barrett}}, \bibinfo {author} {\bibfnamefont {H.}~\bibnamefont {Garcke}}, \
  and\ \bibinfo {author} {\bibfnamefont {R.}~\bibnamefont {Nurnberg}},\
  }\bibfield  {title} {\enquote {\bibinfo {title} {Numerical computations of
  the dynamics of fluidic membranes and vesicles},}\ }\href {\doibase
  10.1103/PhysRevE.92.052704} {\bibfield  {journal} {\bibinfo  {journal} {Phys.
  Rev. E}\ }\textbf {\bibinfo {volume} {92}},\ \bibinfo {pages} {052704}
  (\bibinfo {year} {2015})}\BibitemShut {NoStop}%
\bibitem [{\citenamefont {Canham}(1970)}]{Can-1970}%
  \BibitemOpen
  \bibfield  {author} {\bibinfo {author} {\bibfnamefont {P.}~\bibnamefont
  {Canham}},\ }\bibfield  {title} {\enquote {\bibinfo {title} {The minimum
  energy of bending as a possible explanation of the biconcave shape of the
  human red blood cell},}\ }\href {\doibase 10.1016/S0006-3495(73)86044-8}
  {\bibfield  {journal} {\bibinfo  {journal} {J. Theor. Biol.}\ }\textbf
  {\bibinfo {volume} {26}},\ \bibinfo {pages} {61} (\bibinfo {year}
  {1970})}\BibitemShut {NoStop}%
\bibitem [{\citenamefont {W.Helfrich}(1973)}]{Hel-1973}%
  \BibitemOpen
  \bibfield  {author} {\bibinfo {author} {\bibnamefont {W.Helfrich}},\
  }\bibfield  {title} {\enquote {\bibinfo {title} {Elastic properties of lipid
  bilayers: theory and possible experiments},}\ }\href {\doibase
  10.1515/znc-1973-11-1209} {\bibfield  {journal} {\bibinfo  {journal} {Z.
  Naturforsch. C}\ }\textbf {\bibinfo {volume} {28}},\ \bibinfo {pages} {693}
  (\bibinfo {year} {1973})}\BibitemShut {NoStop}%
\bibitem [{\citenamefont {Evans}(1974)}]{Evans1974923}%
  \BibitemOpen
  \bibfield  {author} {\bibinfo {author} {\bibfnamefont {E.~A.}\ \bibnamefont
  {Evans}},\ }\bibfield  {title} {\enquote {\bibinfo {title} {Bending
  resistance and chemically induced moments in membrane bilayers},}\ }\href
  {\doibase 10.1016/S0006-3495(74)85959-X} {\bibfield  {journal} {\bibinfo
  {journal} {Biophys. J.}\ }\textbf {\bibinfo {volume} {14}},\ \bibinfo {pages}
  {923} (\bibinfo {year} {1974})}\BibitemShut {NoStop}%
\bibitem [{\citenamefont {Laadhari}\ \emph {et~al.}(2010)\citenamefont
  {Laadhari}, \citenamefont {Misbah},\ and\ \citenamefont
  {Saramito}}]{laadhari10}%
  \BibitemOpen
  \bibfield  {author} {\bibinfo {author} {\bibfnamefont {A.}~\bibnamefont
  {Laadhari}}, \bibinfo {author} {\bibfnamefont {C.}~\bibnamefont {Misbah}}, \
  and\ \bibinfo {author} {\bibfnamefont {P.}~\bibnamefont {Saramito}},\
  }\bibfield  {title} {\enquote {\bibinfo {title} {On the equilibrium equation
  for a generalized biological membrane energy by using a shape optimization
  approach},}\ }\href {\doibase 10.1016/j.physd.2010.04.001} {\bibfield
  {journal} {\bibinfo  {journal} {Physica D}\ }\textbf {\bibinfo {volume}
  {239}},\ \bibinfo {pages} {1567} (\bibinfo {year} {2010})}\BibitemShut
  {NoStop}%
\bibitem [{\citenamefont {Cottet}\ \emph {et~al.}(2008)\citenamefont {Cottet},
  \citenamefont {Maitre},\ and\ \citenamefont {Milcent}}]{CotMaiMil2008}%
  \BibitemOpen
  \bibfield  {author} {\bibinfo {author} {\bibfnamefont {G.-H.}\ \bibnamefont
  {Cottet}}, \bibinfo {author} {\bibfnamefont {E.}~\bibnamefont {Maitre}}, \
  and\ \bibinfo {author} {\bibfnamefont {T.}~\bibnamefont {Milcent}},\
  }\bibfield  {title} {\enquote {\bibinfo {title} {Eulerian formulation and
  {L}evel-{S}et models for incompressible fluid-structure interaction},}\
  }\href {\doibase 10.1051/m2an:2008013} {\bibfield  {journal} {\bibinfo
  {journal} {Math. Model. Numer. Anal.}\ }\textbf {\bibinfo {volume} {42}},\
  \bibinfo {pages} {471} (\bibinfo {year} {2008})}\BibitemShut {NoStop}%
\bibitem [{\citenamefont {Laadhari}\ \emph {et~al.}(2014)\citenamefont
  {Laadhari}, \citenamefont {Saramito},\ and\ \citenamefont
  {Misbah}}]{LSM2014}%
  \BibitemOpen
  \bibfield  {author} {\bibinfo {author} {\bibfnamefont {A.}~\bibnamefont
  {Laadhari}}, \bibinfo {author} {\bibfnamefont {P.}~\bibnamefont {Saramito}},
  \ and\ \bibinfo {author} {\bibfnamefont {C.}~\bibnamefont {Misbah}},\
  }\bibfield  {title} {\enquote {\bibinfo {title} {Computing the dynamics of
  biomembranes by combining conservative level set and adaptive finite element
  methods},}\ }\href {\doibase 10.1016/j.jcp.2013.12.032} {\bibfield  {journal}
  {\bibinfo  {journal} {J. Comput. Phys.}\ }\textbf {\bibinfo {volume} {263}},\
  \bibinfo {pages} {328} (\bibinfo {year} {2014})}\BibitemShut {NoStop}%
\bibitem [{\citenamefont {Doyeux}\ \emph {et~al.}(2013)\citenamefont {Doyeux},
  \citenamefont {Guyot}, \citenamefont {Chabannes}, \citenamefont
  {Prud'homme},\ and\ \citenamefont {Ismail}}]{Doyeux2013251}%
  \BibitemOpen
  \bibfield  {author} {\bibinfo {author} {\bibfnamefont {V.}~\bibnamefont
  {Doyeux}}, \bibinfo {author} {\bibfnamefont {Y.}~\bibnamefont {Guyot}},
  \bibinfo {author} {\bibfnamefont {V.}~\bibnamefont {Chabannes}}, \bibinfo
  {author} {\bibfnamefont {C.}~\bibnamefont {Prud'homme}}, \ and\ \bibinfo
  {author} {\bibfnamefont {M.}~\bibnamefont {Ismail}},\ }\bibfield  {title}
  {\enquote {\bibinfo {title} {Simulation of two-fluid flows using a finite
  element/level set method. {A}pplication to bubbles and vesicle dynamics},}\
  }\href {\doibase 10.1016/j.cam.2012.05.004} {\bibfield  {journal} {\bibinfo
  {journal} {J. Comput. Appl. Math.}\ }\textbf {\bibinfo {volume} {246}},\
  \bibinfo {pages} {251} (\bibinfo {year} {2013})}\BibitemShut {NoStop}%
\bibitem [{\citenamefont {Du}\ \emph {et~al.}(2004)\citenamefont {Du},
  \citenamefont {Liu},\ and\ \citenamefont {Wang}}]{DuLiuWan-2004}%
  \BibitemOpen
  \bibfield  {author} {\bibinfo {author} {\bibfnamefont {Q.}~\bibnamefont
  {Du}}, \bibinfo {author} {\bibfnamefont {C.}~\bibnamefont {Liu}}, \ and\
  \bibinfo {author} {\bibfnamefont {X.}~\bibnamefont {Wang}},\ }\bibfield
  {title} {\enquote {\bibinfo {title} {A phase field approach in the numerical
  study of the elastic bending energy for vesicle membranes},}\ }\href
  {\doibase 10.1016/j.jcp.2004.01.029} {\bibfield  {journal} {\bibinfo
  {journal} {J. Comput. Phys.}\ }\textbf {\bibinfo {volume} {198}},\ \bibinfo
  {pages} {450} (\bibinfo {year} {2004})}\BibitemShut {NoStop}%
\bibitem [{\citenamefont {Hu}\ \emph {et~al.}(2014)\citenamefont {Hu},
  \citenamefont {Kim},\ and\ \citenamefont {Lai}}]{Hu2014670}%
  \BibitemOpen
  \bibfield  {author} {\bibinfo {author} {\bibfnamefont {W.-F.}\ \bibnamefont
  {Hu}}, \bibinfo {author} {\bibfnamefont {K.~Y.}\ \bibnamefont {Kim}}, \ and\
  \bibinfo {author} {\bibfnamefont {M.-C.}\ \bibnamefont {Lai}},\ }\bibfield
  {title} {\enquote {\bibinfo {title} {An immersed boundary method for
  simulating the dynamics of three-dimensional axisymmetric vesicles in
  {N}avier-{S}tokes flows},}\ }\href {\doibase 10.1016/j.jcp.2010.03.020}
  {\bibfield  {journal} {\bibinfo  {journal} {J. Comput. Phys.}\ }\textbf
  {\bibinfo {volume} {257, Part A}},\ \bibinfo {pages} {670} (\bibinfo {year}
  {2014})}\BibitemShut {NoStop}%
\bibitem [{\citenamefont {Rahimian}\ \emph {et~al.}(2010)\citenamefont
  {Rahimian}, \citenamefont {Veerapaneni},\ and\ \citenamefont
  {Biros}}]{RahVeeBir-2010}%
  \BibitemOpen
  \bibfield  {author} {\bibinfo {author} {\bibfnamefont {A.}~\bibnamefont
  {Rahimian}}, \bibinfo {author} {\bibfnamefont {S.~K.}\ \bibnamefont
  {Veerapaneni}}, \ and\ \bibinfo {author} {\bibfnamefont {G.}~\bibnamefont
  {Biros}},\ }\bibfield  {title} {\enquote {\bibinfo {title} {Dynamic
  simulation of locally inextensible vesicles suspended in an arbitrary
  two-dimensional domain, a boundary integral method},}\ }\href {\doibase
  10.1016/j.jcp.2010.05.006} {\bibfield  {journal} {\bibinfo  {journal} {J.
  Comput. Phys.}\ }\textbf {\bibinfo {volume} {229}},\ \bibinfo {pages} {6466}
  (\bibinfo {year} {2010})}\BibitemShut {NoStop}%
\bibitem [{\citenamefont {Kaoui}\ \emph {et~al.}(2011)\citenamefont {Kaoui},
  \citenamefont {Harting},\ and\ \citenamefont {Misbah}}]{Kaoui2011}%
  \BibitemOpen
  \bibfield  {author} {\bibinfo {author} {\bibfnamefont {B.}~\bibnamefont
  {Kaoui}}, \bibinfo {author} {\bibfnamefont {J.}~\bibnamefont {Harting}}, \
  and\ \bibinfo {author} {\bibfnamefont {C.}~\bibnamefont {Misbah}},\
  }\bibfield  {title} {\enquote {\bibinfo {title} {Two-dimensional vesicle
  dynamics under shear flow: Effect of confinement},}\ }\href {\doibase
  10.1103/PhysRevE.83.066319} {\bibfield  {journal} {\bibinfo  {journal} {Phys.
  Rev. E}\ }\textbf {\bibinfo {volume} {83}},\ \bibinfo {pages} {066319}
  (\bibinfo {year} {2011})}\BibitemShut {NoStop}%
\bibitem [{\citenamefont {Salac}\ and\ \citenamefont
  {Miksis}(2012)}]{Salac2012}%
  \BibitemOpen
  \bibfield  {author} {\bibinfo {author} {\bibfnamefont {D.}~\bibnamefont
  {Salac}}\ and\ \bibinfo {author} {\bibfnamefont {M.}~\bibnamefont {Miksis}},\
  }\bibfield  {title} {\enquote {\bibinfo {title} {Reynolds number effects on
  lipid vesicles},}\ }\href {\doibase 10.1017/jfm.2012.380} {\bibfield
  {journal} {\bibinfo  {journal} {J. Fluid Mech.}\ }\textbf {\bibinfo {volume}
  {711}},\ \bibinfo {pages} {122} (\bibinfo {year} {2012})}\BibitemShut
  {NoStop}%
\bibitem [{\citenamefont {Laadhari}\ \emph {et~al.}(2017)\citenamefont
  {Laadhari}, \citenamefont {Saramito}, \citenamefont {Misbah},\ and\
  \citenamefont {Szekely}}]{LSMS2018}%
  \BibitemOpen
  \bibfield  {author} {\bibinfo {author} {\bibfnamefont {A.}~\bibnamefont
  {Laadhari}}, \bibinfo {author} {\bibfnamefont {P.}~\bibnamefont {Saramito}},
  \bibinfo {author} {\bibfnamefont {C.}~\bibnamefont {Misbah}}, \ and\ \bibinfo
  {author} {\bibfnamefont {G.}~\bibnamefont {Szekely}},\ }\bibfield  {title}
  {\enquote {\bibinfo {title} {Fully implicit methodology for the dynamics of
  biomembranes and capillary interfaces by combining the level set and newton
  methods},}\ }\href {\doibase 10.1016/j.jcp.2017.04.019} {\bibfield  {journal}
  {\bibinfo  {journal} {J. Comput. Phys.}\ }\textbf {\bibinfo {volume} {343}},\
  \bibinfo {pages} {271} (\bibinfo {year} {2017})}\BibitemShut {NoStop}%
\bibitem [{\citenamefont {Valizadeh}\ and\ \citenamefont
  {Rabczuk}(2022)}]{VALIZADEH2022114191}%
  \BibitemOpen
  \bibfield  {author} {\bibinfo {author} {\bibfnamefont {N.}~\bibnamefont
  {Valizadeh}}\ and\ \bibinfo {author} {\bibfnamefont {T.}~\bibnamefont
  {Rabczuk}},\ }\bibfield  {title} {\enquote {\bibinfo {title} {Isogeometric
  analysis of hydrodynamics of vesicles using a monolithic phase-field
  approach},}\ }\href {\doibase 10.1016/j.cma.2021.114191} {\bibfield
  {journal} {\bibinfo  {journal} {Comput. Methods Appl. Mech. Eng.}\ }\textbf
  {\bibinfo {volume} {388}},\ \bibinfo {pages} {114191} (\bibinfo {year}
  {2022})}\BibitemShut {NoStop}%
\bibitem [{\citenamefont {Laadhari}(2018{\natexlab{a}})}]{LAADHARI2018376}%
  \BibitemOpen
  \bibfield  {author} {\bibinfo {author} {\bibfnamefont {A.}~\bibnamefont
  {Laadhari}},\ }\bibfield  {title} {\enquote {\bibinfo {title} {Implicit
  finite element methodology for the numerical modeling of incompressible
  two-fluid flows with moving hyperelastic interface},}\ }\href {\doibase
  10.1016/j.amc.2018.03.074} {\bibfield  {journal} {\bibinfo  {journal} {Appl.
  Math. Comput.}\ }\textbf {\bibinfo {volume} {333}},\ \bibinfo {pages} {376}
  (\bibinfo {year} {2018}{\natexlab{a}})}\BibitemShut {NoStop}%
\bibitem [{\citenamefont {Laadhari}(2018{\natexlab{b}})}]{LAADHARI201835}%
  \BibitemOpen
  \bibfield  {author} {\bibinfo {author} {\bibfnamefont {A.}~\bibnamefont
  {Laadhari}},\ }\bibfield  {title} {\enquote {\bibinfo {title} {An operator
  splitting strategy for fluid–structure interaction problems with thin
  elastic structures in an incompressible newtonian flow},}\ }\href {\doibase
  10.1016/j.aml.2018.01.001} {\bibfield  {journal} {\bibinfo  {journal} {Appl.
  Math. Lett.}\ }\textbf {\bibinfo {volume} {81}},\ \bibinfo {pages} {35}
  (\bibinfo {year} {2018}{\natexlab{b}})}\BibitemShut {NoStop}%
\bibitem [{\citenamefont {Torres-Sánchez}\ \emph {et~al.}(2019)\citenamefont
  {Torres-Sánchez}, \citenamefont {Millan},\ and\ \citenamefont
  {Arroyo}}]{torres2019}%
  \BibitemOpen
  \bibfield  {author} {\bibinfo {author} {\bibfnamefont {A.}~\bibnamefont
  {Torres-Sánchez}}, \bibinfo {author} {\bibfnamefont {D.}~\bibnamefont
  {Millan}}, \ and\ \bibinfo {author} {\bibfnamefont {M.}~\bibnamefont
  {Arroyo}},\ }\bibfield  {title} {\enquote {\bibinfo {title} {Modelling fluid
  deformable surfaces with an emphasis on biological interfaces},}\ }\href
  {\doibase 10.1017/jfm.2019.341} {\bibfield  {journal} {\bibinfo  {journal}
  {J. Fluid Mech.}\ }\textbf {\bibinfo {volume} {872}},\ \bibinfo {pages}
  {218–271} (\bibinfo {year} {2019})}\BibitemShut {NoStop}%
\bibitem [{\citenamefont {Cokelet}\ \emph {et~al.}(1963)\citenamefont
  {Cokelet}, \citenamefont {Merrill}, \citenamefont {Gilliland}, \citenamefont
  {Shin}, \citenamefont {Britten},\ and\ \citenamefont
  {Jr}}]{cokelet1963rheology}%
  \BibitemOpen
  \bibfield  {author} {\bibinfo {author} {\bibfnamefont {G.~R.}\ \bibnamefont
  {Cokelet}}, \bibinfo {author} {\bibfnamefont {E.}~\bibnamefont {Merrill}},
  \bibinfo {author} {\bibfnamefont {E.}~\bibnamefont {Gilliland}}, \bibinfo
  {author} {\bibfnamefont {H.}~\bibnamefont {Shin}}, \bibinfo {author}
  {\bibfnamefont {A.}~\bibnamefont {Britten}}, \ and\ \bibinfo {author}
  {\bibfnamefont {R.~W.}\ \bibnamefont {Jr}},\ }\bibfield  {title} {\enquote
  {\bibinfo {title} {The rheology of human blood—measurement near and at zero
  shear rate},}\ }\href@noop {} {\bibfield  {journal} {\bibinfo  {journal}
  {Transactions of the Society of Rheology}\ }\textbf {\bibinfo {volume} {7}},\
  \bibinfo {pages} {303} (\bibinfo {year} {1963})}\BibitemShut {NoStop}%
\bibitem [{\citenamefont {Neofytou}\ and\ \citenamefont
  {Drikakis}(2003)}]{NEOFYTOU2003127}%
  \BibitemOpen
  \bibfield  {author} {\bibinfo {author} {\bibfnamefont {P.}~\bibnamefont
  {Neofytou}}\ and\ \bibinfo {author} {\bibfnamefont {D.}~\bibnamefont
  {Drikakis}},\ }\bibfield  {title} {\enquote {\bibinfo {title} {Non-newtonian
  flow instability in a channel with a sudden expansion},}\ }\href {\doibase
  10.1016/S0377-0257(03)00041-7} {\bibfield  {journal} {\bibinfo  {journal} {J.
  Non-Newton. Fluid Mech.}\ }\textbf {\bibinfo {volume} {111}},\ \bibinfo
  {pages} {127} (\bibinfo {year} {2003})}\BibitemShut {NoStop}%
\bibitem [{\citenamefont {Copley}(1960)}]{copley1960apparent}%
  \BibitemOpen
  \bibfield  {author} {\bibinfo {author} {\bibfnamefont {A.}~\bibnamefont
  {Copley}},\ }\href@noop {} {\enquote {\bibinfo {title} {Apparent viscosity
  and wall adherence of blood systems},}\ } (\bibinfo {year}
  {1960})\BibitemShut {NoStop}%
\bibitem [{\citenamefont {{W}illmore}(1993)}]{Wil-1993}%
  \BibitemOpen
  \bibfield  {author} {\bibinfo {author} {\bibfnamefont {T.}~\bibnamefont
  {{W}illmore}},\ }\href@noop {} {\emph {\bibinfo {title} {Riemannian
  geometry}}}\ (\bibinfo  {publisher} {Oxford University Press, USA},\ \bibinfo
  {year} {1993})\ \bibinfo {note} {{ISBN} = 0198514921}\BibitemShut {NoStop}%
\bibitem [{\citenamefont {Feng}\ and\ \citenamefont
  {Klug}(2006)}]{feng-klug-06}%
  \BibitemOpen
  \bibfield  {author} {\bibinfo {author} {\bibfnamefont {F.}~\bibnamefont
  {Feng}}\ and\ \bibinfo {author} {\bibfnamefont {W.~S.}\ \bibnamefont
  {Klug}},\ }\bibfield  {title} {\enquote {\bibinfo {title} {Finite element
  modeling of lipid bilayer membranes},}\ }\href {\doibase
  10.1016/j.jcp.2006.05.023} {\bibfield  {journal} {\bibinfo  {journal} {J.
  Comput. Phys.}\ }\textbf {\bibinfo {volume} {220}},\ \bibinfo {pages} {394}
  (\bibinfo {year} {2006})}\BibitemShut {NoStop}%
\bibitem [{\citenamefont {Laadhari}\ and\ \citenamefont
  {Szekely}(2017)}]{Laadhari20171047}%
  \BibitemOpen
  \bibfield  {author} {\bibinfo {author} {\bibfnamefont {A.}~\bibnamefont
  {Laadhari}}\ and\ \bibinfo {author} {\bibfnamefont {G.}~\bibnamefont
  {Szekely}},\ }\bibfield  {title} {\enquote {\bibinfo {title} {Fully implicit
  finite element method for the modeling of free surface flows with surface
  tension effect},}\ }\href {\doibase 10.1002/nme.5493} {\bibfield  {journal}
  {\bibinfo  {journal} {Int. J. Numer. Methods Eng.}\ }\textbf {\bibinfo
  {volume} {111}},\ \bibinfo {pages} {1047} (\bibinfo {year}
  {2017})}\BibitemShut {NoStop}%
\bibitem [{\citenamefont {Walburn}\ and\ \citenamefont
  {Schneck}(1976)}]{Walburn1976}%
  \BibitemOpen
  \bibfield  {author} {\bibinfo {author} {\bibfnamefont {F.~J.}\ \bibnamefont
  {Walburn}}\ and\ \bibinfo {author} {\bibfnamefont {D.~J.}\ \bibnamefont
  {Schneck}},\ }\bibfield  {title} {\enquote {\bibinfo {title} {A constitutive
  equation for whole human blood},}\ }\href {\doibase 10.3233/bir-1976-13307}
  {\bibfield  {journal} {\bibinfo  {journal} {Biorheology}\ }\textbf {\bibinfo
  {volume} {13}},\ \bibinfo {pages} {201} (\bibinfo {year} {1976})}\BibitemShut
  {NoStop}%
\bibitem [{\citenamefont {{J. Janela}}\ \emph {et~al.}(2005)\citenamefont {{J.
  Janela}}, \citenamefont {{A. Lefebvre}},\ and\ \citenamefont {{B.
  Maury}}}]{janela}%
  \BibitemOpen
  \bibfield  {author} {\bibinfo {author} {\bibnamefont {{J. Janela}}}, \bibinfo
  {author} {\bibnamefont {{A. Lefebvre}}}, \ and\ \bibinfo {author}
  {\bibnamefont {{B. Maury}}},\ }\bibfield  {title} {\enquote {\bibinfo {title}
  {A penalty method for the simulation of fluid - rigid body interaction},}\
  }\href {\doibase 10.1051/proc:2005010} {\bibfield  {journal} {\bibinfo
  {journal} {ESAIM: Proc.}\ }\textbf {\bibinfo {volume} {14}},\ \bibinfo
  {pages} {115} (\bibinfo {year} {2005})}\BibitemShut {NoStop}%
\bibitem [{\citenamefont {Alnaes}\ \emph {et~al.}(2015)\citenamefont {Alnaes},
  \citenamefont {Blechta}, \citenamefont {Hake}, \citenamefont {Johansson},
  \citenamefont {Kehlet}, \citenamefont {A.~Logg}, \citenamefont {Ring},
  \citenamefont {Rognes},\ and\ \citenamefont {Wells}}]{AlnaesEtal2015}%
  \BibitemOpen
  \bibfield  {author} {\bibinfo {author} {\bibfnamefont {M.~S.}\ \bibnamefont
  {Alnaes}}, \bibinfo {author} {\bibfnamefont {J.}~\bibnamefont {Blechta}},
  \bibinfo {author} {\bibfnamefont {J.}~\bibnamefont {Hake}}, \bibinfo {author}
  {\bibfnamefont {A.}~\bibnamefont {Johansson}}, \bibinfo {author}
  {\bibfnamefont {B.}~\bibnamefont {Kehlet}}, \bibinfo {author} {\bibfnamefont
  {C.~R.}\ \bibnamefont {A.~Logg}}, \bibinfo {author} {\bibfnamefont
  {J.}~\bibnamefont {Ring}}, \bibinfo {author} {\bibfnamefont {M.~E.}\
  \bibnamefont {Rognes}}, \ and\ \bibinfo {author} {\bibfnamefont {G.~N.}\
  \bibnamefont {Wells}},\ }\bibfield  {title} {\enquote {\bibinfo {title} {The
  {FEniCS} project version 1.5},}\ }\href {\doibase
  10.11588/ans.2015.100.20553} {\bibfield  {journal} {\bibinfo  {journal}
  {Archive of Numerical Software}\ }\textbf {\bibinfo {volume} {3}} (\bibinfo
  {year} {2015}),\ 10.11588/ans.2015.100.20553}\BibitemShut {NoStop}%
\bibitem [{\citenamefont {LeVeque}(1996)}]{Leveque-1996}%
  \BibitemOpen
  \bibfield  {author} {\bibinfo {author} {\bibfnamefont {R.~J.}\ \bibnamefont
  {LeVeque}},\ }\bibfield  {title} {\enquote {\bibinfo {title} {High-resolution
  conservative algorithms for advection in incompressible flow},}\ }\href
  {\doibase 10.1137/0733033} {\bibfield  {journal} {\bibinfo  {journal} {SIAM
  J. Numer. Anal.}\ }\textbf {\bibinfo {volume} {33}},\ \bibinfo {pages} {627}
  (\bibinfo {year} {1996})}\BibitemShut {NoStop}%
\bibitem [{\citenamefont {Zhao}\ and\ \citenamefont
  {Shaqfeh}(2011)}]{Zhao2011}%
  \BibitemOpen
  \bibfield  {author} {\bibinfo {author} {\bibfnamefont {H.}~\bibnamefont
  {Zhao}}\ and\ \bibinfo {author} {\bibfnamefont {E.~S.~G.}\ \bibnamefont
  {Shaqfeh}},\ }\bibfield  {title} {\enquote {\bibinfo {title} {The dynamics of
  a vesicle in simple shear flow},}\ }\href {\doibase
  10.1017/S0022112011000115} {\bibfield  {journal} {\bibinfo  {journal} {J.
  Fluid Mech.}\ }\textbf {\bibinfo {volume} {674}},\ \bibinfo {pages} {578}
  (\bibinfo {year} {2011})}\BibitemShut {NoStop}%
\bibitem [{\citenamefont {Kraus}\ \emph {et~al.}(1996)\citenamefont {Kraus},
  \citenamefont {Wintz}, \citenamefont {Seifert},\ and\ \citenamefont
  {Lipowsky}}]{Kraus1996}%
  \BibitemOpen
  \bibfield  {author} {\bibinfo {author} {\bibfnamefont {M.}~\bibnamefont
  {Kraus}}, \bibinfo {author} {\bibfnamefont {W.}~\bibnamefont {Wintz}},
  \bibinfo {author} {\bibfnamefont {U.}~\bibnamefont {Seifert}}, \ and\
  \bibinfo {author} {\bibfnamefont {R.}~\bibnamefont {Lipowsky}},\ }\bibfield
  {title} {\enquote {\bibinfo {title} {Fluid vesicles in shear flow},}\ }\href
  {\doibase 10.1103/PhysRevE.72.011901} {\bibfield  {journal} {\bibinfo
  {journal} {Phys. Rev. Lett.}\ }\textbf {\bibinfo {volume} {77}},\ \bibinfo
  {pages} {3685} (\bibinfo {year} {1996})}\BibitemShut {NoStop}%
\bibitem [{\citenamefont {Kantsler}\ and\ \citenamefont
  {Steinberg}(2005)}]{Kantsler2005}%
  \BibitemOpen
  \bibfield  {author} {\bibinfo {author} {\bibfnamefont {V.}~\bibnamefont
  {Kantsler}}\ and\ \bibinfo {author} {\bibfnamefont {V.}~\bibnamefont
  {Steinberg}},\ }\bibfield  {title} {\enquote {\bibinfo {title} {Orientation
  and dynamics of a vesicle in tank-treading motion in shear},}\ }\href
  {\doibase 10.1103/PhysRevLett.95.258101} {\bibfield  {journal} {\bibinfo
  {journal} {Phys. Rev. Lett.}\ }\textbf {\bibinfo {volume} {95}},\ \bibinfo
  {pages} {258101} (\bibinfo {year} {2005})}\BibitemShut {NoStop}%
\bibitem [{\citenamefont {Fischer}\ \emph {et~al.}(1978)\citenamefont
  {Fischer}, \citenamefont {St\"ohr-Liesen},\ and\ \citenamefont
  {Schmid-Sch\"onbein}}]{Fischer-1978}%
  \BibitemOpen
  \bibfield  {author} {\bibinfo {author} {\bibfnamefont {T.~M.}\ \bibnamefont
  {Fischer}}, \bibinfo {author} {\bibfnamefont {M.}~\bibnamefont
  {St\"ohr-Liesen}}, \ and\ \bibinfo {author} {\bibfnamefont {H.}~\bibnamefont
  {Schmid-Sch\"onbein}},\ }\bibfield  {title} {\enquote {\bibinfo {title} {The
  red cell as a fluid droplet: Tank tread-like motion of the human erythrocyte
  membrane in shear flow},}\ }\href@noop {} {\bibfield  {journal} {\bibinfo
  {journal} {Science}\ }\textbf {\bibinfo {volume} {202}},\ \bibinfo {pages}
  {894} (\bibinfo {year} {1978})}\BibitemShut {NoStop}%
\bibitem [{\citenamefont {Gizzi}\ \emph {et~al.}(2015)\citenamefont {Gizzi},
  \citenamefont {Ruiz-Baier}, \citenamefont {Rossi}, \citenamefont {Laadhari},
  \citenamefont {Cherubini},\ and\ \citenamefont {Filippi}}]{Gizzi2015157}%
  \BibitemOpen
  \bibfield  {author} {\bibinfo {author} {\bibfnamefont {A.}~\bibnamefont
  {Gizzi}}, \bibinfo {author} {\bibfnamefont {R.}~\bibnamefont {Ruiz-Baier}},
  \bibinfo {author} {\bibfnamefont {S.}~\bibnamefont {Rossi}}, \bibinfo
  {author} {\bibfnamefont {A.}~\bibnamefont {Laadhari}}, \bibinfo {author}
  {\bibfnamefont {C.}~\bibnamefont {Cherubini}}, \ and\ \bibinfo {author}
  {\bibfnamefont {S.}~\bibnamefont {Filippi}},\ }\bibfield  {title} {\enquote
  {\bibinfo {title} {A three-dimensional continuum model of active contraction
  in single cardiomyocytes},}\ }\href {\doibase 10.1007/978-3-319-05230-4_6}
  {\bibfield  {journal} {\bibinfo  {journal} {Modeling, Simulation and
  Applications}\ }\textbf {\bibinfo {volume} {14}},\ \bibinfo {pages} {157}
  (\bibinfo {year} {2015})}\BibitemShut {NoStop}%
\end{thebibliography}%

\end{document}